\documentclass[a4paper,fleqn]{cas-sc}

\usepackage[numbers]{natbib}
\usepackage{hyperref}
\usepackage{amsmath, amsthm, amsfonts}
\usepackage{tikz}
\usepackage[font=small,labelfont=bf]{caption}

\usepackage{float} 
\usepackage{booktabs}
\usepackage{microtype} 
\usepackage{siunitx}
\usepackage{verbatim}
\usepackage{xcolor}

\usepackage{etoolbox}
\apptocmd{\normalsize}{%
  \setlength{\abovedisplayskip}{3pt}%
  \setlength{\belowdisplayskip}{3pt}%
  \setlength{\abovedisplayshortskip}{0pt}%
  \setlength{\belowdisplayshortskip}{0pt}%
}

\DeclareSIUnit\angstrom{\text{\AA}} 

\newcommand{\diff}{\mathop{}\!\mathrm{d}}

\def\tsc#1{\csdef{#1}{\textsc{\lowercase{#1}}\xspace}}
\tsc{MD}
\tsc{YS}
\tsc{EA}
\tsc{SZ}

\begin{document}
\let\WriteBookmarks\relax
\def\floatpagepagefraction{0.8}
\def\textpagefraction{.001}

\shorttitle{Polar solvation energy of nonrigid biomolecules}    

\shortauthors{M. Dogbatsey et al.}  

\title [mode = title]{On calculating polar solvation energy of nonrigid proteins in the Poisson-Boltzmann theory}

%

\author[1]{Matthias Dogbatsey}[orcid = 0000-0001-6905-6684]





\affiliation[1]{organization={Department of Mathematics, University of Alabama},
            city={Tuscaloosa},
            postcode={35487}, 
            state={AL},
            country={USA}}

\author[1]{Yuanzhen Shao}[orcid = 0000-0002-5348-7153]





\author[2,3,4]{Emil Alexov}[orcid = 0000-0001-5346-0156]

\affiliation[2]{organization={Department of Physics and Astronomy, Clemson University},
            city={Clemson},
            postcode={29634}, 
            state={SC},
            country={USA}}

\affiliation[3]{organization={Medical Biophysics Graduate Program, Clemson University},
            addressline={118 Kinard Laboratory}, 
            city={Clemson},
            postcode={29634}, 
            state={SC},
            country={USA}}

\affiliation[4]{organization={Clemson University Institute for Human Genetics},
            addressline={106 Gregor Mendel Circle}, 
            city={Greenwood},
            postcode={29646}, 
            state={SC},
            country={USA}}

\author[1]{Shan Zhao}[orcid = 0000-0002-3023-2107]


\cormark[1]

\ead{szhao@ua.edu}




\cortext[1]{Corresponding author}



\begin{abstract}
The Poisson-Boltzmann (PB) theory is a cornerstone of implicit solvent models for electrostatic analysis, and has found a great success in various biomolecular applications. However, in calculating polar solvation energy, 
one should consider that the structure of the protein changes upon transition from vacuum to water phases. To address this, here we report for the first time 
a generalized PB framework capable of accommodating nonrigid conformational changes without suffering from self-energy artifacts. 
For regularized PB models, in which the charge singularities are captured by the Green's functions, self-energies in the water and vacuum states will be analytically canceled. For non-regularized PB solvers, such as APBS and DelPhi, a simple thermodynamic cycle is proposed for nonrigid proteins by adding a Coulombic correction in vacuum. The generalized PB theory is validated using a perturbed two-atom system and a diverse set of proteins with different structures in vacuum and water, demonstrating its accuracy and robustness, regardless of the choice of sharp-interface and diffuse-interface PB models and different numerical solvers.
\end{abstract}



\maketitle

\section{Introduction}

The Poisson-Boltzmann (PB) model is one of the most widely used implicit solvent models for electrostatic analysis, especially for calculating the polar solvation energy, or the electrostatic component of the free energy of placing biomolecules, e.g. proteins, from vacuum to water phases. 
However, even though proteins are not rigid, their structures in the PB solvation energy calculation are assumed to be fixed in solvent and vacuum states. This is because the PB energy functional involves self-energy terms \cite{baker2005improving,zhao2024calculation}, i.e., the singular charge at an atom center will interact with the potential induced by itself, which artificially yields infinite energies. Consequently, a rigid structure has to be assumed in the PB theory, so that the same self-energies involved in solvent and vacuum states can be canceled.  
The rigidity assumption of the PB theory may be partially bypassed if the polar solvation energy can be computed without using singular potential values near charge centers, such as in the induced surface charges method \cite{rocchia2002rapid} for the sharp-interface PB model and in the far-field energy method \cite{ijaodoro2025generalizing} for heterogeneous PB models. Nevertheless, to the best of our knowledge, the direct treatment of nonrigid structures in the PB solvation energy calculation has never been tackled. 

In this work, an analytical treatment of nonrigid structures will be first formulated for the regularized PB models \cite{geng2017two,zhao2024calculation}. Moreover, a simple thermodynamic cycle will be proposed to model nonrigid proteins by adding a correction term. This enables the commonly used non-regularized PB solvers, such as APBS \cite{jurrus2018improvements} and DelPhi \cite{li2012delphi}, to handle nonrigid proteins in polar solvation energy calculation, through post-processing. For simplicity, the proposed developments will be presented for a diffuse-interface PB model \cite{zhao2024calculation}, and can be similarly applied to the sharp-interface PB model \cite{geng2017two}.

\section{Analytical treatment of nonrigid proteins in a regularized PB model}
Consider a protein comprising $N_m$ atoms. 
When immersed in an ionic solvent, for the $j$th atom with a radius $R_j$, a partial charge with a valence $z_j$ is assumed to be located at the atom center $\mathbf{r}_j$, $j=1, \ldots, N_m$. 
This study assumes that the protein structure in the vacuum is different, i.e., the atom center becomes $\hat{\mathbf{r}}_j$, with $z_j$ and $R_j$ unchanged. 
We denote the solute structure in the solvent and vacuum by $\mathcal{M} = \{\mathbf{r}_j, z_j\}_{j=1}^{N_m} $ and $\hat{\mathcal{M}} = \{\hat{\mathbf{r}}_j, z_j\}_{j=1}^{N_m}$, respectively.

The PB energy functionals in terms of the dimensionless potential for the diffuse interface PB (DIPB) model have been presented in \cite{zhao2024calculation}. Following their notation, consider a diffuse interface characterized by a surface function $S(\mathbf{r})$, generated by the Gaussian convolution surface (GCS) algorithm \cite{zhao2024calculation}. The surface function $S(\mathbf{r})$ takes constant values $S=1$ and $S=0$, respectively, in the solute region $\Omega_i$ and solvent region $\Omega_e$. In the smooth solute-solvent boundary region $\Omega_t$, $S(\mathbf{r})$ varies smoothly from 1 to 0. The smooth dielectric function can then be calculated as $\epsilon(S) = \epsilon_m S + \epsilon_s (1-S)$ over the entire domain $\Omega$, where $\epsilon_m = 1$ and $\epsilon_s = 80$. 

The dimensionless potential $u$ in the solvent state is governed by the nonlinear PB (NPB) equation,
\begin{equation}
-\nabla \cdot (\epsilon \nabla u) + (1-S)\kappa^2 \sinh(u) = \rho :=4\pi \frac{e_c^2}{k_B T} \sum_{j=1}^{N_m} z_j \delta(\mathbf{r} -\mathbf{r}_j),  \label{eq:PBE}
\end{equation}
subject to a Dirichlet boundary condition. 
Here $\kappa$, $k_B$, $T$ and $e_c$ denote the modified Debye-H\"uckel parameter, Boltzmann constant, absolute temperature, and fundamental charge, respectively. In the vacuum state, the dimensionless potential $v$ satisfies a Poisson equation 
\begin{equation}
-\epsilon_m \Delta v = \hat{\rho} : = 4\pi \frac{e_c^2}{k_B T} \sum_{j=1}^{N_m} z_j \delta(\mathbf{r} -\hat{\mathbf{r}}_j),\label{eq:PE}
\end{equation}
with a Dirichlet boundary condition. 
We note that the sharp-interface PB model can also be formulated in terms of Eq. \eqref{eq:PBE}, by defining $S(\mathbf{r})$ as a Heaviside function  \cite{shao2023convergence}.

The polar solvation energy of transferring the protein from vacuum to the solvent environment is the difference between the electrostatic free energies in the solvent and vacuum states. Based on the energy functionals given in \cite{zhao2024calculation}, the PB solvation energy for perturbed molecular structures $\mathcal{M}$ and $\hat{\mathcal{M}}$ can be defined as
\begin{align}
    \Delta E &= \dfrac{1}{2} k_B T \sum_{j=1}^{N_m} z_j (u(\mathbf{r}_j) - v(\hat{\mathbf{r}}_j)) + I_{u}(u) \label{DE-DIPB}\\
    I_{u}(u) &= - \dfrac{1}{4\pi} \dfrac{(k_BT)^2}{e_c^2} \int_{\Omega_i^c} (1-S) \kappa^2 (\cosh(u) - 1) \, \diff\mathbf{r} + \dfrac{1}{8\pi} \dfrac{(k_BT)^2}{e_c^2} \int_{\Omega_i^c} (1-S) \kappa^2 u \sinh(u) \,\diff\mathbf{r}, \nonumber
\end{align}
where $\Omega_i^c$ is the complement of the solute domain $\Omega_i$. By using a trilinear approximation of singular sources $\rho$ and $\hat{\rho}$, the finite difference method can be applied to solve potentials $u$ and $v$. However, due to the mismatch of charge centers in the first term of \eqref{DE-DIPB},  self-energies associated with solvent and vacuum states cannot be canceled. Consequently, $\Delta E$ calculated by the DIPB model will diverge.


The state-of-the-art numerical treatment of the PB singular charges is using regularization schemes \cite{chern2003accurate,geng2017two,lee2021regularization,wang2021regularization}, in which $u$ is decomposed into two or three parts so that the singular component can be analytically captured by the Green's function. However, no one has attempted to rigorously treat nonrigid proteins in the PB regularization before. 

The regularization scheme for the diffuse interface PB model was first developed in \cite{wang2021regularization}.
Following that, the potential in the water state is decomposed as the sum of the reaction field and Coulombic potentials, i.e., $u = u_{RF} + u_C$. The singular part $u_C$ satisfies a Poisson equation with the singular source $\rho$ and can be solved analytically as the Green's function:
$u_C (\mathbf{r}) =  G (\mathbf{r} )= \frac{e_c^2}{k_B T} \sum_{j=1}^{N_m} \frac{z_j}{\epsilon_m |\mathbf{r} - \mathbf{r}_j|}$. The other component $u_{RF}$ satisfies a regularized NPB equation 
\cite{wang2021regularization},
\begin{equation}
-\nabla \cdot (\epsilon \nabla u_{RF}) + (1-S)\kappa^2 \sinh(u_{RF} + G) = \nabla \epsilon \cdot \nabla G, \label{eq:PBE_REG}
\end{equation}
with a smooth and bounded source. 
In the vacuum state, the potential $v$ due to the perturbed singular source $\hat{\rho}$ can also be represented analytically as the Green's function $v(\mathbf{r}) =  \hat{G} (\mathbf{r} )= \frac{e_c^2}{k_B T} \sum_{j=1}^{N_m} \frac{z_j}{\epsilon_m |\mathbf{r} - \hat{\mathbf{r}}_j |}$.
Then, the polar solvation energy for perturbed molecular structures $\mathcal{M}$ and $\hat{\mathcal{M}}$ can be given as
\begin{align}
\Delta E &= \frac{1}{2} k_B T \int_{\Omega} \sum_{i=1}^{N_m} z_i \delta(\mathbf{r} -\mathbf{r}_i) (u_{RF}(\mathbf{r}) + G(\mathbf{r})) \, \diff\mathbf{r} - \frac{1}{2} k_B T \int_{\Omega} \sum_{i=1}^{N_m} z_i \delta(\mathbf{r} -\hat{\mathbf{r}}_i) \Hat{G}(\mathbf{r}) \, \diff\mathbf{r} + I_u(u_{RF}+G)\nonumber  \\
&= \frac{1}{2} k_B T \sum_{i=1}^{N_m} z_i u_{RF}(\mathbf{r}_i) + I_u(u_{RF}+G) 
+ \frac{1}{2} e_c^2\sum_{i=1}^{N_m} \sum_{\substack{j=1 }}^{N_m} \frac{z_i z_j}{\epsilon_m |\mathbf{r}_i - \mathbf{r}_j|} - \frac{1}{2} e_c^2\sum_{i=1}^{N_m} \sum_{\substack{j=1 }}^{N_m} \frac{z_i z_j}{\epsilon_m |\hat{\mathbf{r}}_i - \hat{\mathbf{r}}_j|}, \label{DE-REG-0}
\end{align}
where the definitions of Green's functions $G$ and $\hat{G}$ have been applied. The first two terms of Eq. \eqref{DE-REG-0} are \emph{bounded}, while the singular self-energies are obviously presented in the last two terms when $j=i$. 

In this study, we propose a regularization (REG) scheme to calculate $\Delta E$. We note that the self-energies of the water and vacuum states can actually be canceled in Eq. \eqref{DE-REG-0}, even though the atom centers $\mathbf{r}_i$ are different from $\hat{\mathbf{r}}_i$. Thus, $\Delta E$ can be calculated as
\begin{equation}
\Delta E = \frac{1}{2} k_B T \sum_{i=1}^{N_m} z_i u_{RF}(\mathbf{r}_i) + I_u(u_{RF}+G) 
+ \frac{1}{2} e_c^2\sum_{i=1}^{N_m} \sum_{\substack{j=1 \\ j \neq i}}^{N_m} \frac{z_i z_j}{\epsilon_m |\mathbf{r}_i - \mathbf{r}_j|} - \frac{1}{2} e_c^2\sum_{i=1}^{N_m} \sum_{\substack{j=1 \\ j \neq i}}^{N_m} \frac{z_i z_j}{\epsilon_m |\hat{\mathbf{r}}_i - \hat{\mathbf{r}}_j|}.
\label{DE-REG}
\end{equation}
We note that for calculating the first two terms of Eq. \eqref{DE-REG}, one just needs to solve the regularized NPB equation for $u_{RF}$ in the water state, similar to the regularization scheme for the rigid proteins \cite{wang2021regularization}. After removing self-energies, the last two terms of Eq. \eqref{DE-REG} are actually the Coulombic energies of $\mathcal{M}$ and $\hat{\mathcal{M}}$, respectively, and they can be computed analytically. 

\begin{figure}[pos=t]
    \centering
    \includegraphics[width=0.6\linewidth]{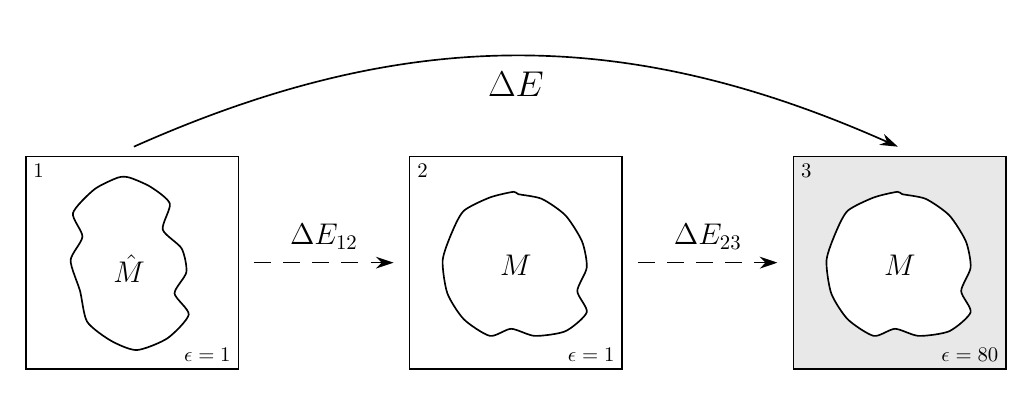}
\caption{A thermodynamic cycle is proposed for calculating the polar solvation energy $\Delta E$ for an object that adopts different shapes in vacuum and water phases. The white background indicates the vacuum. The gray background indicates a medium with a dielectric constant $\epsilon_s = 80$. Here $\Delta E = \Delta E_{12} + \Delta E_{23}$, with 
$\Delta E_{23}$ being the polar solvation energy of a rigid biomolecular structure, and $\Delta E_{12}$ being the Coulombic energy required to induce the conformational change in the vacuum.}
    \label{fig:conformational_change}
\end{figure}

\section{Correction for non-regularized PB solvers}
Motivated by Eq. \eqref{DE-REG} of the REG scheme, we propose a simple thermodynamic cycle to model $\Delta E$ for non-regularized PB solvers, see Fig. \ref{fig:conformational_change}. An intermediate State 2 is inserted in between the vacuum State 1 and the solvent State 3, with the vacuum environment but with the molecular structure $\mathcal{M}$ as in the solvent.
The protein is assumed to have conformational changes in the vacuum between State 1 and State 2. The corresponding polar energy is actually the difference of Coulombic energies in the two states and is given analytically by
\begin{align}
    \Delta E_{12} =  \frac{1}{2} e_c^2\sum_{i=1}^{N_m} \sum_{\substack{j=1 \\ j \neq i}}^{N_m} \frac{z_i z_j}{\epsilon_m |\mathbf{r}_i - \mathbf{r}_j|} - \frac{1}{2} e_c^2\sum_{i=1}^{N_m} \sum_{\substack{j=1 \\ j \neq i}}^{N_m} \frac{z_i z_j}{\epsilon_m |\hat{\mathbf{r}}_i - \hat{\mathbf{r}}_j|}. \label{eq:DeltaE12}
\end{align}
For the solvation process between State 2 and State 3, a rigid molecular structure can be assumed, so that any existing PB solver can be employed to approximate $\Delta E_{23}$. Then, the polar solvation energy for perturbed molecular structures $\mathcal{M}$ and $\hat{\mathcal{M}}$ is calculated as $\Delta E = \Delta E_{12} + \Delta E_{23}$. Here, $\Delta E_{12}$ can be regarded as a correction, which can be applied to any PB solver as a post-processing. 

For the corrected diffuse-interface PB (C-DIPB) scheme, one first solves $u$ in the water by the NPB equation \eqref{eq:PBE}. In the vacuum, one solves a Poisson equation with the same structures $\mathcal{M}$ for $v$
\begin{equation}
-\epsilon_m \Delta v = 4\pi \frac{e_c^2}{k_B T} \sum_{j=1}^{N_m} z_j \delta(\mathbf{r} -\mathbf{r}_j). \label{eq:PE_OG}
\end{equation}
As $u$ and $v$ have singularities at the same atom centers $\mathbf{r}_j$, the self-energies can be canceled in 
\begin{align*}
\Delta E_{23} &= \frac{1}{2} k_B T \sum_{j=1}^{N_m} z_j (u(\mathbf{r}_j) - v(\mathbf{r}_j)) + I_u(u). 
\end{align*}
This correction procedure can be adapted to any PB solver, including regularized or non-regularized and sharp-interface or diffuse-interface ones. In this study, the correction to APBS \cite{jurrus2018improvements} and DelPhi \cite{li2012delphi} packages will be considered, and the resulted schemes will be called C-APBS and C-DelPhi, respectively. We note that the APBS package directly accepts different molecular structures for water and vacuum states. Such results will be referred to as APBS energies. 

\section{Numerical verification of the effectiveness of the correction procedure}
Six PB schemes will be examined in this section. Based on the diffuse-interface PB model, the DIPB and REG schemes are coded in one package, by using central finite differences over a uniform mesh with an isotropic grid spacing $h$ in $x$, $y$, and $z$ directions. The modified Debye-H\"uckel  parameter $\kappa$ is computed as $\kappa^2 \approx 8.487\, \text{\AA}^{-2} I$, where $I$ represents the ionic strength of the solvent. Throughout our simulations, we set the solvent dielectric $\epsilon_s = 80$, the solute dielectric $\epsilon_m = 1$, the solvent probe radius to $1.4\,\text{\AA}$, and the ionic strength to $I = 0.15\,\text{M}$. Other setups such as Dirichlet boundary conditions are similarly handled as in Ref. \cite{zhao2024calculation,ijaodoro2025generalizing}. The C-DIPB result is simply obtained by adding a correction $\Delta E_{12}$ given by Eq. \eqref{eq:DeltaE12}. 
The DelPhi, C-DelPhi, APBS and C-APBS schemes are based on the sharp-interface NPB equation, with model parameters being similarly chosen. Solvation energies were converted using standard thermodynamic equivalents ($1\,\text{kcal/mol} = 4.184\,\text{kJ/mol}$ for APBS, and $1\,\text{kT} = 0.5922\,\text{kcal/mol}$ at $298\,\text{K}$ for DelPhi). 

\textbf{Example 1:} Consider a two-atom system with equal radii $R_1 = R_2 = \qty{2}{\angstrom}$ and $(z_1,z_2) =(-1,1)$. In the water, the atom centers are fixed at $\mathbf{r}_1 = (-2.2, 0, 0)$ and $\mathbf{r}_2 = (2.2, 0, 0)$, while $\hat{\mathbf{r}}_1$ and $\hat{\mathbf{r}}_2$ are randomly perturbed in the vacuum
\begin{equation}
    \hat{\mathbf{r}}_j = \mathbf{r}_{j} - \frac{C_p}{2} (1,1,1) + C_p (\gamma_{j,1}, \gamma_{j,2} , \gamma_{j,3}), 
 \quad \text{for } j = 1, 2, \label{eq:Perturbed}
\end{equation}
where $C_p$ is the perturbation length, and $\gamma_{j,i}$ is drawn from a uniform distribution $\mathcal{U}(0,1)$. The perturbation in Eq. \eqref{eq:Perturbed} is shifted such that the positions $\hat{\mathbf{r}}_j$ are centered around the original coordinates $\mathbf{r}_{j}$. A fixed computational domain $\Omega=[-8,8]^3$ is used in all computations. 

\begin{table}[htbp]
\centering
\caption{The NPB polar solvation energy for a two-atom system with perturbed atom centers. The results by varying $C_p$ for a fixed $h=0.6$ and varying $h$ for a fixed $C_p=0.1$ are listed.}
\label{tab:perturbation}
\begin{tabular}{l r r r r r | l r r r r r}
\toprule
\multicolumn{6}{c|}{$h=0.6$} & \multicolumn{6}{c}{$C_p=0.1$}  \\
\cline{1-6} \cline{7-12}
$C_p$ & DIPB & APBS & REG & C-DIPB & C-APBS & $h$ & DIPB & APBS & REG & C-DIPB & C-APBS\\ 
\hline
0.1 & -79.86  & -101.25 & -66.48 &-66.06 &  -90.45 & 0.8 & -60.73  & -112.02 & -52.83 & -51.60 &  -101.14 \\ 
0.2 & -79.12  & -15.89  & -65.38 &-64.96 &  -89.35 & 0.6 & -79.86  & -101.25 & -66.48 & -66.06 &  -90.45 \\ 
0.3 & -86.58  & 218.50  & -65.29 &-64.87 &  -89.26 & 0.4 & 233.70  & -74.21  & -68.82 & -68.28 &  -89.78 \\ 
0.4 & -176.29 & 141.55  & -67.40 &-66.98 &  -91.37 & 0.2 & 1971.66 & 513.99  & -72.42 & -72.70 &  -89.42 \\ 
0.5 & -96.00  & 207.37  & -61.05 &-60.63 &  -85.02 & 0.1 & 5715.48 & 2109.99 & -74.47 & -74.94 &  -86.72 \\  
\bottomrule
\end{tabular}%
\end{table}

The NPB solvation energies are reported in Table \ref{tab:perturbation}. It is clear that the uncorrected DIPB and APBS methods fail, because the PB self-energies cannot be canceled once the structural rigidity assumption is relaxed. In particular, when $h$ becomes smaller with $C_p=0.1$, the self-energies become larger and larger, so that the DIPB and APBS energies diverge drastically. When $h$ is fixed to be a large number, such as $h=0.6$, the self-energies of DIPB and APBS are not very large. Nevertheless, the DIPB and APBS energies are sensitively affected when $C_p$ becomes larger. 
After corrections, both C-DIPB and C-APBS energies are bounded as $h$ goes to zero. When $C_p$ increases with $h=0.6$, both C-DIPB and C-APBS energies change robustly. Note that the C-DIPB energies are very close to those of REG in all cases, indicating a high accuracy of the correction procedure.

\textbf{Example 2:} We next study a set of proteins with different structures in the water and vacuum states. The set of 74 proteins studied in \cite{chakravorty2018reproducing} is initially chosen. The NPB equation is solved in all cases, and four proteins (PDB ID: 1MC2, 2FDN, 2H5C, and 4TKB) are excluded in this study because some NPB solvers fail to converge for them. For the remaining set of 70 proteins, their vacuum structures $\hat{\mathcal{M}}$ are assumed to be their crystal structures from the protein databank. The molecular structures in the water state $\mathcal{M}$ are taken as the energy-minimized structures in the presence of explicit TIP3P water molecules \cite{chakravorty2018reproducing}.
For the DIPB, C-DIPB, and REG, we set a probe radius of \qty{1.5}{\angstrom}, $h=\qty{0.5}{\angstrom}$ and an edge value of \qty{5}{\angstrom} \cite{zhao2024calculation}. 
For the DelPhi calculations, we applied a grid scale of $2.0\,\text{grids}/\text{\AA}$ and a percentage fill (\texttt{perfil}) of $70$~\cite{chakravorty2018reproducing}. 
For APBS, computations were executed by using either a $97 \times 97 \times 97$ grid or a grid spacing of $0.3 \times 0.3 \times 0.3\,\text{\AA}$. 
Further details on the setup used in APBS and DelPhi packages and the correction postprocessing in energy calculation can be found on GitHub. 

\begin{table}[!th]
\centering
\caption{The NPB polar solvation energies for 10 nonrigid proteins. For each protein, $\hat{\mathcal{M}}$ in the vacuum is assumed to be the crystal structure, while $\mathcal{M}$ in the water is the energy-minimized structure in the presence of explicit TIP3P water molecules \cite{chakravorty2018reproducing}.}
\label{tab:TIP3P_XTAL}
\begin{tabular}{l r r r r r r}
\toprule
PDB ID & DIPB &  APBS &  REG & C-DIPB & C-APBS & C-DelPhi \\ 
\midrule
1TG0 & -227.58 & -3685.02 & -2458.15 & -2490.82 & -2341.84 & -2705.76 \\ 
1TQG & -904.50 & 166.80 & -1244.60 & -1289.37 & -712.18 & -1602.61 \\ 
1VBW & -2466.56 & -1104.90 & -1321.86 & -1352.82 & -1353.26 & -1544.71 \\ 
1W0N & 263.74 & -671.87 & -1385.22 & -1424.97 & -1207.14 & -1728.95 \\ 
1IQZ & -3674.82 & -8851.22 & -3835.97 & -3871.99 & -3532.39 & -4110.90 \\ 
3LZT & 557.40 & 1317.13 & -1431.82 & -1476.58 & -1283.77 & -1781.26 \\ 
3O5Q & -2399.62 & -1142.01 & -1682.90 & -1736.62 & -1229.46 & -2113.36 \\ 
3PUC & -43.56 & 1171.87 & -954.36 & -990.04 & -777.72 & -1273.57 \\ 
3VOR & -3280.79 & -5629.64 & -1284.89 & -1337.06 & -833.55 & -1793.57 \\ 
4A02 & 265.52 & -1356.83 & -957.31 & -1004.40 & -598.23 & -1375.92 \\ 
\bottomrule
\end{tabular}%
\end{table}

\begin{figure}[pos=ht]
    \centering
    \includegraphics[width=0.8\linewidth]{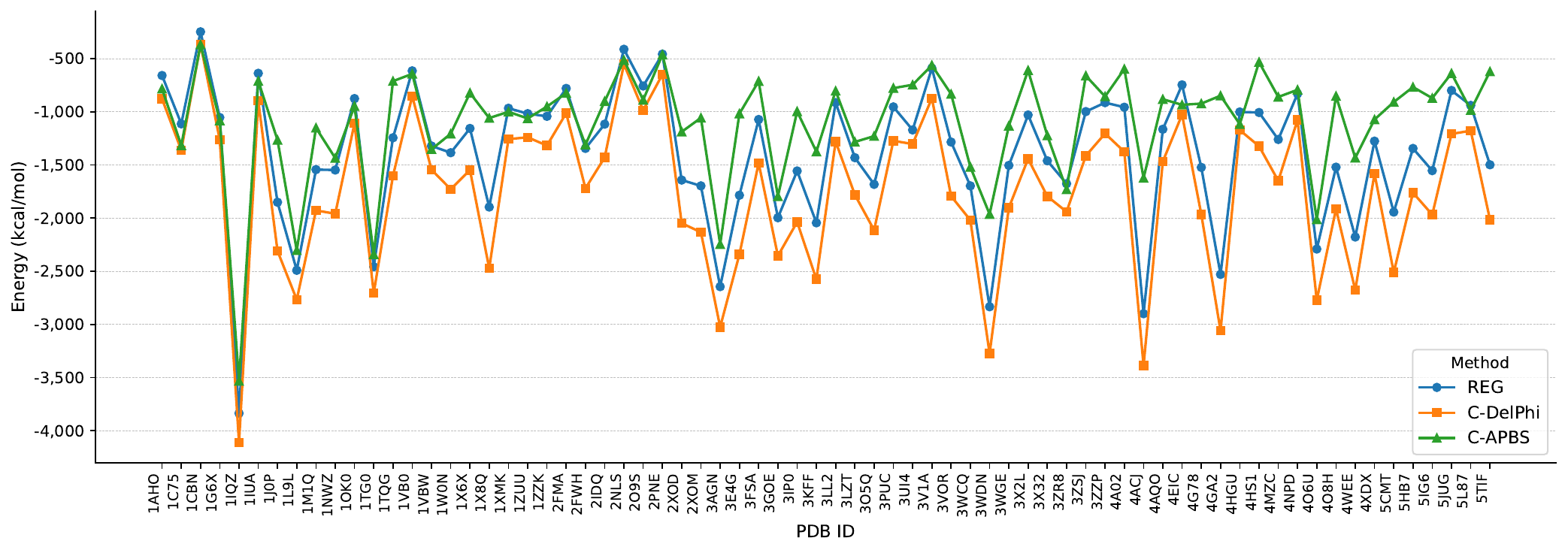}
\caption{The NPB solvation energies of 70 perturbed proteins. The Pearson correlation coefficients are $0.99$ and $0.83$, respectively, for C-DelPhi vs REG and C-APBS vs REG (with $p$-values of $3.29 \times 10^{-61}$ and $3.42 \times 10^{-19}$).}
    \label{fig:protein_energy}
\end{figure}

Out of the set of 70 proteins, 10 are selected to report in Table~\ref{tab:TIP3P_XTAL} for their NPB solvation energies generated by six schemes. It is obvious that the uncorrected dual-structure evaluations in DIPB and APBS exhibit profound deviations from others, and some even become unphysical with positive values. After correction, the C-DIPB energies are quite close to the REG energies, because they are based on the same diffuse-interface PB model. Based on sharp-interface PB model, C-APBS and C-DelPhi energies are different, due to different numerical solvers. But both packages successfully deliver reasonable energies after correction. To further analyze C-APBS and C-DelPhi results, their NPB solvation energies of 70 proteins are plotted against those of REG in Fig. \ref{fig:protein_energy}. 
It can be seen that the correlation between REG and C-DelPhi energies is stronger than that between REG and C-APBS energies. Indeed, the Pearson correlation coefficients are found to be $0.99$ and $0.83$, respectively, for  C-DelPhi vs REG and C-APBS vs REG, with vanishing $p$-values. 
This suggests that the C-DelPhi provides more accurate results than the C-APBS. 
Overall, the high correlations in Fig. \ref{fig:protein_energy}  demonstrate the success of the proposed correction, regardless of the underlying PB models and numerical solvers. 


\section{Conclusion}

In this study, we introduced a generalized PB theory to accommodate conformational changes in calculating polar solvation energy. When the charge singularities are represented through Green's functions in the regularization, the PB self-energies in the water and vacuum states can be analytically canceled. For non-regularized PB models, a simple thermodynamic cycle is proposed by introducing a structural perturbation in vacuum, while a rigid structure can still be assumed in the solvation calculation. This generalized energy correction is versatile and can be readily integrated into standard PB solvers like APBS and DelPhi. Numerical validations on a perturbed two-atom system and a comprehensive set of $70$ proteins demonstrate the effectiveness of the generalized PB theory based on different PB models and numerical solvers. Moreover, the C-DelPhi is shown to be more accurate than the C-APBS. 

\bigskip
\noindent
{\bf Acknowledgments}\\
This research was partially supported by the National Science Foundation (NSF) under grants DMS-2512104, DMS-2512105, and DMS-2306991, and by the National Institutes of Health (NIH) under grant R35GM151964.

\bigskip
\noindent
{\bf Data availability}\\
The protein structures analyzed in this study were downloaded from the protein databank https://www.rcsb.org, and processed according to the cited references. The correction procedure used in C-APBS and C-DelPhi energy calculation is available on GitHub: https://github.com/szhao-ua/PB-solvation-energy-of-nonrigid-proteins. 

\bibliographystyle{cas-model2-names}

\bibliography{nonrigid}

\end{document}